\documentclass[11pt]{article}
\usepackage{amsxtra,amssymb,amsthm,amsmath,latexsym}
\def\oH{\buildrel\circ\over H}
\def\oH1{\buildrel\circ\over H\kern-.02in{}^1}

\def\ep{\epsilon}

\def\d{{\delta}}

\begin{document}


\title{
Dynamical Systems Method for ill-posed equations with monotone 
operators
\footnote{Math subject classification: 34R30,  35R25, 35R30,
37C35, 37L05, 37N30, 47A52,  
47J06, 65M30, 65N21$\quad$ key words: dynamical systems method, 
ill-posed problems, monotone operators, iterative methods }
}

\author{ A.G. Ramm\\
Mathematics Department, Kansas State University, \\
Manhattan, KS 66506-2602, USA\\
E:mail: ramm@math.ksu.edu
}

\date{}
\maketitle

\begin{abstract}
 Consider an operator equation (*) $B(u)-f=0$ in a real Hilbert space.
 Let us call this equation ill-posed if the  operator
 $B'(u)$ is not boundedly invertible, and well-posed otherwise.
The DSM (dynamical systems method) for solving equation (*)
consists of a construction of a Cauchy 
problem, which has the following properties:
 1) it has a global solution for an arbitrary initial data,
 2) this solution tends to a limit as time tends to infinity,
 3) the limit is the minimal-norm solution to the equation
$B(u)=f$.\\ \hspace*{0.5cm} A global convergence theorem is
proved for DSM for equation (*) with monotone $C_{loc}^2$
operators $B$. \end{abstract}




\section{Introduction}
In this paper the dynamical systems method, DSM, (see [1]) 
is proposed for
solving nonlinear operator
equation of the form:
$$
B(u)-f=0, \quad f\in H.
\eqno{(1.1)}
$$
We make the following assumptions: 

A)\,\,{\it $B$ is a monotone, nonlinear, $C^2_{loc}$
operator in a real Hilbert space $H$, i.e, $\sup_{u\in
{\mathcal B}(u_0,R)}||B^{(j)}(u)||\leq M_j(R):=M_j,\,j=0, 1,
2,$ where $R>0$ is arbitrary, ${\mathcal B}(u_0,R):=\{u:
||u-u_0||\leq R \}$, $B^{(j)}(u)$ is the Fr\'echet
derivative. The set
 $N:=\{z: B(z)-f=0\}$ is non-empty, and $y\in N$ 
is its minimal-norm 
element. 
}

It is known
that $N$ is convex and closed if $B$ is monotone and hemicontinuous, and 
such  $N\in H$ 
contains the unique minimal-norm 
element $y$: $B(y)-f=0, \quad ||y||\leq ||z||,\,\, \forall z\in N.$

Let $\dot u$ denote the derivative with respect to time. Consider the 
dynamical system ( that is, the Cauchy problem ):
$$
\dot u=\Phi(u), \,\,\, u(0)=u_0;\,\,\quad \Phi:=-A_\ep^{-1}[B(u)+\ep u-f], 
\eqno{(1.2)}
$$
where $u:=u_\ep(t)$,
 $A_\ep:=A+\ep I,\, A:=B'(u)$, $I$ is the identity operator,
$\ep>0$ is a number,  the prime stands for the  Fr\'echet
derivative, and
$\Phi(u)$ is locally Lipschitz if the assumption A) holds.
Thus, problem (1.2) has a unique local solution.

{\it The DSM in this paper consists of solving equation
(1.1) by solving (1.2), and proving that for any initial
approximation $u_0$ the following results (1.3) and (1.4)
hold:} $$ \exists u(t) \forall t>0, \quad \exists
V_\ep:=u(\infty):=\lim_{t\to \infty}u(t), \quad B(V_\ep)+\ep
V_\ep-f=0, \eqno{(1.3)} $$ and $$\lim_{\ep\to
0}||V_\ep-y||=0. \eqno{(1.4)} $$ Conclusion (1.4) is known,
but we give a simple proof in Section 2 for convenience of
the reader. Denote $$ t_\ep:=-2\log(\ep). \eqno{(1.5)} $$
{\it We prove:} $$ \lim_{\ep \to 0}||u_\ep(t_\ep)-y||=0.
\eqno{(1.6)} $$ Problem (1.1) with noisy data $f_\delta$,
$||f_\delta -f||\leq \delta$, given in place of $f$,
generates the problem: $$ \dot w_\d=\Phi_\d(w_\d), \,\,\,
w_\d(0)=u_0, \eqno{(1.7)} $$ and $\Phi_\d$ is similar to
$\Phi$ in (1.2) with $f_\d$ and $w_\d$ replacing $f$ and
$u$, respectively. We prove that, for a suitable choice of
$\ep(\delta)$ and $t_\delta$, the solution $w_\d$ to (1.7),
calculated at a suitable stopping time
$t=t_\d:=t_{\ep(\d)}$, converges to $y$: $$ \lim_{\d \to
0}||w_\d(t_\d)-y||=0, \eqno{(1.8)} $$ where $\ep(\d)>0, \,
\lim_{\d\to 0}\ep(\d)=0,$ and the choice of $t_\d$ with
property (1.8) is called the stopping rule. {\it Thus,
calculating the solution to (1.7) at a suitable stopping
time $t_\delta$ yields a stable approximation to the
minimal-norm solution $y$ of equation (1.1) in the sense
(1.8).}

{\it We do not restrict the growth of nonlinearity at
infinity and do not assume that the initial approximation
$u_0$ is close to the solution $y$ in any sense.}

 Usually convergence theorems for Newton-type methods for
solving nonlinear equation (1.1), even in the well-posed
case, have the assumption that the initial data $u_0$ is
close to $y$. We obtain a {\it global} convergence result
for a continuous regularized Newton-type method (1.2).  In
Section 2 the results are stated and proved, and in Section
3 a convergent iterative process for solving (1.1) is
constructed.  In [1] other versions of DSM are presented.
See also [2].

 \section{Statement and Proof of the Results }

{\bf Theorem 2.1.} {\it If A) holds, then, for any choice of
$u_0$, relations (1.3), (1.4), and (1.6) hold. If $w_\d(t)$
solves (1.7), then there is a $t_\d$ such that (1.8) holds.
}

\medskip

{\bf Proof of Theorem 2.1.}
Let $||B(u)+\ep u-f||:=g(t)$, $u=u_\ep(t)$. Using (1.2), one gets $g\dot g=
-g^2$, so $g(t)=g_0\exp(-t)$, where $g_0:=g(0)$. This and (1.2) imply
$||\dot u||\leq g_0\ep^{-1}\exp(-t)$, because $||A_\ep^{-1}||\leq \ep^{-1}$.
 Therefore $u(\infty):=V_\ep$ does exist by the Cauchy 
criterion, one has:
$$||u(t)-V_\ep||\leq  g_0\ep^{-1}\exp(-t),
\eqno{(2.1)}
$$ 
and the last equation in (1.3) holds. 
Using (1.5) and (2.1), one gets
$$||u(t_\ep)-V_\ep||\leq  g_0\ep.  
\eqno{(2.2)}  
$$
To prove (1.4), use (1.1) and  write equation (1.3) as 
$B(V_\ep)-B(y)+\ep V_\ep=0$, multiply
by $V_\ep-y$, use the monotonicity of $B$, and get
$$
  (V_\ep,V_\ep-y)\leq 0.
\eqno{(2.3)}
$$
Thus, $||V_\ep||\leq ||y||$. Consequently, 
$V_\ep\rightharpoonup v$, as $\ep\to 0$,
where $v$ is
some element and $\rightharpoonup $ stands for weak convergence.
 Since monotone 
hemicontinuous operators are $w-$closed, that is, 
$V_\ep\rightharpoonup v$
and $B(V_\ep)+\ep V_\ep\to f$ imply $B(v)+\ep v=f$, and since the 
inequality
 $||V_\ep||\leq ||y||$ implies $||v||\leq \liminf_{\ep\to 0}||V_\ep||\leq 
||y||,$
one concludes that $v$ is the minimal norm solution 
to (1.1), so $v=y$.
To prove that $ V_\ep\to y$, where $\to$ denotes the strong convergence in 
$H$,
 one uses inequality (2.3) in the form 
$$
(V_\ep-y,V_\ep-y)\leq (y,y-V_\ep),
\eqno{(2.4)}
$$ 
and the weak convergence of $V_\ep$ to $y$. This 
yields statement (1.4).

 Relation (1.6) follows from (1.4) and (2.2). Indeed:
$$ ||u_\ep(t_\ep)-y||\leq ||u_\ep(t_\ep)-V_\ep||+||V_\ep-y||\to 0, \quad \ep\to 0.
$$
Let us prove (1.8). As above (cf (2.1)-(2.2)), one proves
$$
||w_\d(t)-W_\d||\leq  g_{0\d}\ep^{-1}\exp(-t),\,\, 
||w_\d(t_\ep)-W_\d||\leq  g_{0\d}\ep,
\eqno{(2.5)}
$$ 
where $ g_{0\d}:=||B(u_0)+\ep u_0-f_\d||$, and $W_\d:=W_{\d \ep}:=
w_\d(\infty)=\lim_{t\to \infty}w_\d(t)$.
 In general, $f_\d$ is not in 
the range of $B$, and $W_{\d \, \ep}(\infty)$ does not 
converge as $\ep\to 0$, $\d>0$ being fixed.
One can use (2.2) to prove that there is an $\ep:=\ep(\d)\to 0$
as $\d\to 0$, such that (1.8) holds. Indeed, $W_\d:=W_{\d \ep}$ solves the 
equation:
$$
B(W_\d)+\ep W_\d-f_\d=0.
\eqno{(2.6)}
$$
Subtract from (2.6) equation (1.3), denote $W_\d-V_\ep:=\psi_\d:=\psi_{\d \ep},\, 
f_\d-f:=h_\d, \,\, ||h_\d||\leq \d$, multiply the resulting equation by 
$\psi_\d$, and use the monotonicity of $B$ to get $\ep (\psi_\d, 
\psi_\d)\leq \d 
||\psi_\d||$, so 
$$||W_\d-V_\ep||\leq \d \ep^{-1}.
\eqno{(2.7)}
$$ 

Let us assume:
$$
\lim_{\d\to 0}\d \ep^{-1}=0, \quad \lim_{\d\to 0}\ep(\d)=0.
\eqno{(2.8)}
$$
For example, one may take $\ep=\d^b,\,\, 0<b<1$. If (2.8) holds,
then
$$
\lim_{\d\to 0}||W_\d-V_{\ep(\d)}||=0.
  \eqno{(2.9)}
$$
From (2.5),(2.7)-(2.9), and (1.4), one gets (1.8). $\Box$

\medskip

{\bf Remark 2.2.} Using a different argument, the author has 
recently proved (see [2]), that $u_{\ep(t)}(t)\in 
{\mathcal B}(u_0,R),\,\,R:=3(||u_0||+||y||),\,\,
{\mathcal B}(u_0,R):=\{u: ||u-u_0||\leq R\},$ for any initial 
approximation 
$u_0$ and a 
suitable choice of $\ep(t)$.

\medskip

 \section{Convergent Iterative Process}

In this section it is proved that for a suitable choice of
$h_n>0$ and $\ep_n>0, \,\, \lim_{n\to \infty}\ep_n=0,$ the
iterative process: 
$$
u_{n+1}=u_n-h_n[A(u_n)+\ep_n]^{-1}[B(u_n)+\ep_n u_n-f],
\eqno{(3.1)} 
$$ 
where $u_0$ is arbitrary, converges to $y$,
the minimal norm solution of (1.1), provided assumptions A)
from Section 1 hold. In particular, $h_n$ can be a fixed
positive constant, suitably chosen. A possible choice
of $h_n$ and $\ep_n$ is made in the proof of Theorem 3.1 
 
\medskip

{\bf Theorem 3.1.} {\it One can choose $\ep_n$ and $h_n$ so that
$$
\lim_{n\to \infty}||u_n-y||=0,
\eqno{(3.2)} 
$$
where $y$ is the minimal norm solution to (1.1).}

\medskip

{\bf Proof.} Let $V_n:=V_{\ep_n}$ solve the equation:
$$
B(V_n)+\ep_n V_n-f=0,\quad \ep_n>0, \,\, \lim_{n\to \infty}\ep_n=0.
\eqno{(3.3)}
$$
Since $B\in C_{loc}^2$ is monotone, equation (3.3) has a unique solution 
for any 
$f\in H$. Denote $z_n:=u_n-V_n$, and $||z_n||:=g_n$. Then, by the triangle 
inequality, one gets:
$$
||u_n-y||\leq ||u_n-V_n||+||V_n-y||=||z_n||+||V_n-y||.
\eqno{(3.4)}
$$
 In Section 2 we have proved that 
$$
\lim_{n\to \infty}||V_n-y||=0.
\eqno{(3.5)}
$$
Therefore Theorem 3.1 is proved if one proves that:
$$
\lim_{n\to \infty}||z_n||=0.
\eqno{(3.6)}
$$
{\it Let us prove} (3.6). Let
$$
b_n:=||V_{n+1}-V_n||, \quad \lim_{n\to \infty}b_n=0.
\eqno{(3.7)}
$$
 Rewrite (3.1) as
$$
z_{n+1}=(1-h_n)z_n-h_n[A(u_n)+\ep_n]^{-1}K(z_n)- (V_{n+1}-V_n),
\eqno{(3.8)}
$$
where we have used the Taylor formula:
$$
B(u_n)+\ep_n u_n-f=B(u_n)-B(V_n)+\ep_n z_n=[A(u_n)+\ep_n]z_n+K(z_n),
\eqno{(3.9)}
$$
with
$$
||K(z_n)||\leq 0.5 M_2 ||z_n||^2:=cg_n^2, \quad c:=0.5 M_2.
\eqno{(3.10)}
$$
From (3.10) and (3.8) one gets:
$$
g_{n+1}\leq (1-h_n)g_n+ \frac {ch_n} {\ep_n} g_n^2 +b_n, \,\,\, 0<h_n\leq 1.
\eqno{(3.11)}
$$
Let 
$$
\ep_n=2cg_n.
\eqno{(3.12)}
$$
Then (3.11) can be written as:
$$
g_{n+1}\leq (1-a_n)g_n+ b_n, \quad 0<a_n:=0.5 h_n\leq 0.5.
\eqno{(3.13)}
$$
Theorem 3.1 follows now from lemma 3.2:

\medskip

{\bf Lemma 3.2.} {\it If $b_k\geq 0,\, \lim_{k\to
\infty}b_k=0, \,\,
 0<a_n\leq 0.5,$ and $$ \sum_{n=1}^\infty a_n=\infty, \quad
\lim_{n\to \infty}\sum_{k=1}^{n-1} b_k e^{-\sum_{j=k+1}^n
a_j}=0, \eqno{(3.14)} $$ then (3.13) implies $$ \lim_{n\to
\infty} g_n=0. \eqno{(3.15)} $$ }

\medskip

{\bf Proof of lemma 3.2.} From (3.13) by induction one gets:
$$
g_{n+1}\leq b_n+\sum_{k=1}^{n-1}b_k \prod_{j=k+1}^n 
(1-a_j)+g_1\prod_{j=1}^n (1-a_j).
\eqno{(3.16)}
$$
If $a\geq 0$, then $1-a\leq e^{-a}$. Thus, 
$$
 \prod_{j=k+1}^n(1-a_j)\leq e^{-\sum_{j=k+1}^n a_j},\quad 0<a_j<1.
\eqno{(3.17)}
$$
The assumptions of lemma 3.2 and formulas (3.16) and (3.17) imply (3.15). 
$\Box$

{\it In conclusion note that $a_n$, and therefore $h_n=2a_n$, can
always be chosen so that (3.14) hold.}

For example, one can choose  $a_j=\log p_j=\log p, \, 1<p\leq \sqrt{e}$.
This choice corresponds to the choice of constant $h_n:=h:=2\log p$.
Then $ e^{-\sum_{j=k+1}^n a_j}=p^{-n+k}$,
$\sum_{k=1}^{n-1}b_k  e^{-\sum_{j=k+1}^n a_j}=\sum_{k=1}^{n-1}b_k 
p^{-n+k}$, and \\
$\lim_{n\to \infty }\sum_{k=1}^{n-1}b_kp^{-n+k}=0$
provided that $\lim_{k\to \infty }b_k=0$. $\Box$

\end{document}